\documentclass[11pt]{amsart}

\usepackage{amscd,amssymb,amsopn,amsmath,amsthm,mathrsfs,graphics,amsfonts,enumerate,verbatim,calc
}

\usepackage[OT2,OT1]{fontenc}
\newcommand\cyr{%
\renewcommand\rmdefault{wncyr}%
\renewcommand\sfdefault{wncyss}%
\renewcommand\encodingdefault{OT2}%
\normalfont
\selectfont}
\DeclareTextFontCommand{\textcyr}{\cyr} 

\usepackage{amssymb,amsmath}

\DeclareFontFamily{OT1}{rsfs}{}
\DeclareFontShape{OT1}{rsfs}{n}{it}{<-> rsfs10}{}
\DeclareMathAlphabet{\mathscr}{OT1}{rsfs}{n}{it}

\numberwithin{equation}{section}
\newtheorem{theorem}{Theorem}[section]
\newtheorem{lemma}[theorem]{Lemma}
\newtheorem{proposition}[theorem]{Proposition}
\newtheorem{corollary}[theorem]{Corollary}

\newtheorem*{maintheorem}{Main Theorem}
\newtheorem{conjecture}{Conjecture}

\theoremstyle{definition}
\newtheorem*{ack}{Acknowledgement}
\newtheorem{definition}[theorem]{Definition}
\newtheorem{remark}[theorem]{Remark}

\theoremstyle{remark}

\renewcommand{\ker}{\operatorname{Ker}}

\newcommand{\Spec}{\operatorname{Spec}}

\newcommand{\Frac}{\operatorname{Frac}}

\newcommand{\fm}{\frak{m}}

\newcommand{\fp}{\frak{p}}

\newcommand{\Frob}{\operatorname{Frob}}

\begin{document}
\title[An application of the almost purity theorem to the homological conjectures]
{An application of the almost purity theorem to the homological conjectures}

\author[K.Shimomoto]{Kazuma Shimomoto}
\address{Department of Mathematics, College of Humanities and Sciences,
Nihon University, Setagaya-ku, Tokyo, 156-8550, Japan}
\email{shimomotokazuma@gmail.com}
\thanks{The author is partially supported by Grant-in-Aid for Young Scientists (B) \# 25800028}

\subjclass{13A35, 13B22, 13B40, 13D22, 13K05}

\keywords{Almost purity theorem, big Cohen-Macaulay algebra,  Frobenius map, Witt vectors.}


\begin{abstract} 
The aim of this article is to establish the existence of a big Cohen-Macaulay algebra over a local ring in mixed characteristic $p>0$ in the case where the local ring is finite \'etale over a regular local subring after inverting $p$. The main result follows from the so-called almost purity theorem proved by Davis and Kedlaya.
\end{abstract}

\maketitle

\section{Introduction}

The homological conjectures, which consist of a set of statements for commutative Noetherian (local) rings, have been of interest in the quest of commutative algebra (see \cite{Ho07} and \cite{Rob12} for known results and history). In this article, we would like to consider the following conjecture.

\begin{conjecture}
[Hochster]
\label{conjecture}
A Noetherian local ring $(R,\fm,k)$ admits an $R$-algebra $B$ such that $\fm B \ne B$ and every system of parameters of $R$ is a regular sequence on $B$.
\end{conjecture}

We are interested in the situation where the module-finite map of rings $R \to S$ for a regular ring $R$ is \'etale after inverting a prime integer $p>0$. This comes from the condition imposed on the statement of the \textit{almost purity theorem}, which was originally proved by Faltings in $p$-adic Hodge theory, with refinements by Scholze \cite{Sch}, Davis and Kedlaya \cite{DavKed} more recently. A proof given by Davis and Kedlaya uses a valuative method and their original motivation was to pin down the class of rings satisfying the \textit{Witt-perfect} condition.

A very important aspect of the almost purity theorem is expressed by the following facts:
\begin{enumerate}
\item[$\bullet$]
If $A$ is a $p$-torsion free Witt-perfect normal ring and $A \to B$ is an almost \'etale map, then the Frobenius endomorphism on $B/pB$ is also surjective (see Theorem \ref{almostpure}). 

\item[$\bullet$]
If $A$ is an almost Cohen-Macaulay Witt-perfect algebra and $A \to B$ is an almost \'etale map, then $B$ is also almost Cohen-Macaulay (see Proposition \ref{algebramodification}).
\end{enumerate}

We are going to use the version by Davis and Kedlaya, as their statement assumes only $p$-torsion freeness on rings. Our proof reduces to the case where local rings are complete in the topology defined by the powers of the maximal ideals, as is commonly done in the research of commutative algebra. We state our main theorem (see Theorem \ref{maintheorem}).

\begin{maintheorem}
Let $R$ be a regular local ring of mixed characteristic $p>0$ and let $S$ be a torsion free module-finite $R$-algebra such that the localization $R[\frac{1}{p}] \to S[\frac{1}{p}]$ is \'etale. Then $S$ has a balanced big Cohen-Macaulay $R$-algebra. That is, there is an $S$-algebra $T$ such that every system of parameters of $R$ is a regular sequence on $T$.
\end{maintheorem}

In the above setting, $S$ is a semilocal ring. If $S$ is local, then an inspection of the proof of the main theorem shows that every system of parameters of $S$ is a regular sequence on $T$. Assume in general that $R \to S$ is a torsion free module-finite extension of Noetherian domains of generically characteristic 0. Then there is a nonzero element $f \in R$ such that $R[\frac{1}{f}] \to S[\frac{1}{f}]$ is \'etale. Indeed, this follows from the lemma of generic freeness together with the vanishing of the differential module: $\Omega_{\Frac(S)/\Frac(R)}=0$. Thus, our theorem treats the special case that $f=p$ and $R$ is regular.

Faltings' work on $p$-adic Hodge theory led to \textit{almost ring theory}. Its idea turned out to be essential in Heitmann's work on the \textit{Direct Summand Conjecture} in dimension three based on the construction of an \textit{almost Cohen-Macaulay algebra} in \cite{He1}. Later on, Hochster has constructed a big Cohen-Macaulay algebra in dimension three in \cite{Ho02} (see also \cite{AsSh} and \cite{Shim} for relevant works). To prove the main theorem, we employ some tools from \cite{Ho02}. Namely, we first establish the existence of an almost Cohen-Macaulay algebra. Then we prove that it maps to a big Cohen-Macaulay algebra by applying Hochster's \textit{partial algebra modifications}. The following corollary, which is regarded as the special case of the Direct Summand Conjecture, is a consequence of the main theorem (see Corollary \ref{maincorollary}):

\begin{corollary}
\label{DirectSummand}
Let $R$ be a Noetherian regular domain and let $S$ be a torsion free module-finite $R$-algebra. Assume that $R$ is $p$-torsion free and the localization $R[\frac{1}{p}] \to S[\frac{1}{p}]$ is finite \'etale for some prime integer $p>0$. Then $R \hookrightarrow S$ splits as an $R$-module homomorphism.
\end{corollary}

We should point out that Bhatt recently proved Corollary \ref{DirectSummand} in a slightly different setting in \cite{Bha1}. Corollary \ref{DirectSummand} assumes the base ring to be an arbitrary Noetherian regular domain, while Bhatt's result treats a more ramified extension than ours under the assumption that the base ring is \'etale over a ring that is essentially of finite over a discrete valuation ring. More recently, Corollary \ref{DirectSummand} has been generalized in the context of logarithmic algebraic geometry by Gabber and Ramero \cite{GR2}, which is the most general form of the Direct Summand Conjecture proven at this point.

\section{Notation and convention}

All rings in this article are commutative with unity. However, we do not always assume rings to be Noetherian. A \textit{local ring} is a Noetherian ring with a unique maximal ideal. A sequence of elements $(x_1,\ldots,x_n)$ of a local ring $(R,\fm)$ is a \textit{regular sequence} on an $R$-module $M$, if 
\begin{enumerate}
\item[$\bullet$]
$(x_1,\ldots,x_n)M \ne M$.

\item[$\bullet$]
The multiplication map $M/(x_1,\ldots,x_i)M \xrightarrow{x_{i+1}} M/(x_1,\ldots,x_i)M$ is injective for all $i$.
\end{enumerate}

Here is a definition of (balanced) big Cohen-Macaulay algebra or module.

\begin{definition}
Let $M$ be an $R$-algebra (resp. $R$-module) over a local ring $(R,\fm)$. Then we say that $M$ is a \textit{big Cohen-Macaulay R-algebra} (resp. \textit{R-module}), if there is a system of parameters of $R$ that is a regular sequence on $M$. Furthermore, $M$ is \textit{balanced}, if every system of parameters of $R$ is a regular sequence on $M$. 
\end{definition}

There are many flavors of Witt vectors in the literature, however the only Witt vectors we consider in this article are the $p$-typical Witt vectors. The basic reference for $p$-typical Witt vectors is Serre's book \cite{Se}. Another good source is an expository paper \cite{Rab}. We will need to consider the Witt vectors for $p$-torsion free rings to define the Witt-perfect condition. Since the basic part of the Witt vectors is rather involved and we do not require an explicit descriptions, we give only the summary to the extent we need.

Fix a prime integer $p>0$ and a commutative ring $A$ of arbitrary characteristic. Then we denote by $\mathbf{W}(A)$ (resp. $\mathbf{W}_{p^n}(A)$) the ring of (\textit{p-typical}) \textit{Witt vectors} (resp. \textit{Witt vectors of length n}). Then one has a set-theoretic identity: $\mathbf{W}_{p^n}(A)=A^{n+1}$ and a ring-theoretic isomorphism: $\mathbf{W}(A) \cong \varprojlim_n \mathbf{W}_{p^n}(A)$. On the Witt vectors, there is a well-defined ring homomorphism called the \textit{Witt-Frobenius map}: 
$$
\mathbf{F}:\mathbf{W}_{p^n}(A) \to \mathbf{W}_{p^{n-1}}(A),
$$ 
which is described as follows. If $A$ is a ring of prime characteristic $p>0$, then $\mathbf{F}(r_1,r_2,\ldots,r_{n+1})=(r_1^p,r_2^p,\ldots,r_n^p)$ (\cite{DavKed}; Remark 1.5). When $A$ is a general $\mathbf{Z}$-algebra, then the formula is more complicated and found in (\cite{DavKed}; Lemma 1.4). The symbol "Frob" will denote the map $x \mapsto x^p$ for $x \in R$, where $R$ is a ring of characteristic $p>0$ to distinguish it from the similar map $\mathbf{F}$ as above. We note that the set of maps $\mathbf{F}$ with all $n \ge 0$ collate to define a ring homomorphism:
$$
\mathbf{F}:\mathbf{W}(A) \to \mathbf{W}(A),
$$
which we again call the \textit{Witt-Frobenius map}. There are also additive \textit{Verschiebung maps}:
$$
\mathbf{V}:\mathbf{W}_{p^n}(A) \to \mathbf{W}_{p^{n+1}}(A)
$$
defined by $\mathbf{V}(a_0,\ldots,a_n)=(0,a_0,\ldots,a_n)$. We use the notation $\Frac(A)$ for the total ring of fractions for a ring $A$. That is, $\Frac(A)$ is the localization of $A$ with respect to the multiplicative set of nonzero divisors in $A$.

\begin{definition}
Let $p>0$ be a prime number. Say that a ring $A$ has \textit{mixed characteristic $p>0$, or p-torsion free}, if $p$ is a nonzero divisor in $A$.
\end{definition}

\begin{definition}
An $\mathbf{F}_p$-algebra $B$ is \textit{perfect} (resp. \textit{semiperfect}), if the Frobenius map $\Frob:B \to B$ is bijective (resp. surjective). 
\end{definition}

Suppose that $B$ is a perfect $\mathbf{F}_p$-domain which is not a field. Then it is seen that $B$ is not Noetherian.

\begin{definition}[\cite{Ar}]
Let $A$ be an integral domain. The \textit{absolute integral closure} of $A$, denoted by the symbol $A^+$, is defined to be the integral closure of $A$ in a fixed algebraic closure of the field of fractions of $A$.
\end{definition}

\section{Preliminaries for the proof of the main theorem}

First, we construct a $p$-torsion free big ring and use it to construct an almost Cohen-Macaulay algebra whose basic properties will be discussed in the next section. One requires separate treatments in the unramified and ramified cases. Let $k$ be a perfect field of characteristic $p>0$.

\begin{enumerate}
\item[$\bullet$]
Let $(R,\fm)$ be a complete regular local ring of mixed characteristic $p>0$ with residue field $k$ such that $p \notin \fm^2$. Then we can present $R$ as the ring:
$$
\mathbf{W}(k)[[x_2,\ldots,x_d]].
$$
Fix a parameter $\pi$ of the valuation ring $\mathbf{W}(k)$. We define
$$
R_{p^{\infty}}:=\bigcup_{n>0} R_{p^n},
$$
where
$$
R_{p^n}:=\mathbf{W}(k)[\pi_n][[x_2^{p^{-n}},\ldots,x_d^{p^{-n}}]]
$$
and $\{\pi_n\}_{n \in \mathbf{N}}$ and $\{x_i^{p^{-n}}\}_{n \in \mathbf{N}}$ are fixed sequences of elements contained in $R^+$ such that $\pi_0=\pi$ and $\pi_{n+1}^p=\pi_n$.

\item[$\bullet$]
Let $(R,\fm)$ be a complete regular local ring of mixed characteristic $p>0$ with residue field $k$ such that $p \in \fm^2$. Let $T=\mathbf{W}(k)[[t_1,\ldots,t_d]]$. Then we can present $R$ as the quotient ring:
$$
\mathbf{W}(k)[[t_1,\ldots,t_d]]/(p-G)
$$
for some element $G \in \fm_T^2\backslash pT$ for the maximal ideal $\fm_T$ of $T$. We define
$$
R_{p^{\infty}}:=\bigcup_{n>0} R_{p^n},
$$
where
$$
R_{p^n}:=\mathbf{W}(k)[[t_1^{p^{-n}},\ldots,t_d^{p^{-n}}]]/(p-G).
$$

\end{enumerate}

In both unramified and ramified cases, the ring $R_{p^\infty}$ enjoys the following properties (see also Proposition \ref{WittCM}).

\begin{lemma}
\label{WittNormal}
$R_{p^\infty}$ is a non-Noetherian normal domain with a unique maximal ideal and it is contained in $R^+$ and integral, faithfully flat over $R$. Moreover, $\Frob:R_{p^{\infty}}/pR_{p^{\infty}} \to R_{p^{\infty}}/pR_{p^{\infty}}$ is surjective, but not injective. 
\end{lemma}

\begin{proof}
Because $R_{p^\infty}$ is the ascending union of complete regular local rings $R_{p^n}$ with each transition map being module-finite, we see that $R_{p^\infty}$ has a unique maximal ideal and it is integral and faithfully flat over $R$. Moreover, the maximal ideal of $R_{p^{\infty}}$ is not finitely generated. Quite generally, it is known that the filtered colimit of normal domains with injective transition maps is a normal domain (see \cite{SwHu}; Proposition 19.3.1). Thus, $R_{p^\infty}$ is a normal domain. The Frobenius map is surjective on $R_{p^\infty}/pR_{p^\infty}$ by the construction. To see that the Frobenius is not injective, it suffices to prove that $R_{p^\infty}/pR_{p^\infty}$ is not reduced. In the unramified case, consider the sequence $\{\pi_n\}_{n \in \mathbf{N}}$. In the ramified case, consider the sequence $\{G_n\}_{n \in \mathbf{N}}$ by letting $G_n:=G(t_1^{p^{-n}},\ldots,t_d^{p^{-n}})$. Then they form nilpotent elements in $R_{p^\infty}/pR_{p^\infty}$. This proves the lemma.
\end{proof}

This is a typical example of \textit{Witt-perfect ring}, which will be discussed later.

\begin{remark}
A ring map $A \to B$ is \textit{\'etale}, if it is of finite presentation, flat and unramified. Assume that $A \to B$ is \'etale. If $A/pA$ is a (semi)perfect $\mathbf{F}_p$-algebra, then $B/pB$ is also (semi)perfect. This fact holds more generally for an \textit{absolutely flat} map $A \to B$ (see \cite{Oli} for its definition and \cite{GR}; Theorem 3.5.13, or \cite{KedRuo}; Lemma 3.1.5 for a proof). \cite{KedRuo} contains a comprehensive discussion on the tensor category of \'etale algebras over a fixed base ring.
\end{remark}

We recall Hochster's partial algebra modifications. The following discussion is taken from \cite{Ho02}.

Let $R$ be a ring and let $M$ be an $R$-module. Let $M[X_1,\ldots,X_k]:=M \otimes_R R[X_1,\ldots,X_k]$. For an integer $N>0$, let $M[X_1,\ldots,X_k]_{\le N}$ be the $R$-submodule spanned by the elements $uX_1^{a_1}\cdots X_k^{a_k}$ such that $\sum_i a_i \le N$. Now fix an integer $k \ge 0$ and let $x_{k+1}u_{k+1}=\sum_{i=1}^kx_iu_i$ be a relation for $u_i \in M$ and a fixed system of elements $x_1,\ldots,x_k$ (which we will later take to be part of a system of parameters of a local ring). A \textit{partial algebra modification} of $M$ is defined as an $R$-module map:
$$
M \to M':=\frac{M[X_1,\ldots,X_k]_{\le N}}{FM[X_1,\ldots,X_k]_{\le N-1}}
$$
with $F:=(u_{k+1}-\sum_{i=1}^k x_i X_i)$. This process implies that the map $M \to M'$ trivializes the relation $x_{k+1}u_{k+1}=\sum_{i=1}^kx_iu_i$ and it is expected that it will be useful for constructing big Cohen-Macaulay algebras. Let us specialize this construction. Let $(R,\fm)$ be a Noetherian local ring and $x_1,\ldots,x_d$ a system of parameters of $R$ and let $T=M$ be an $R$-algebra. Applying the above construction, we get a sequence of partial algebra modifications:
$$
T \to M_1 \to M_2 \to \cdots \to M_r
$$
We say that the above sequence is \textit{bad}, if $\fm M_r=M_r$. Here is a key lemma for the proof of the main theorem.

\begin{lemma}[\cite{Ho02}; Lemma 5.1]
\label{lemma3}
Let $M$ be a module over a local ring $(R,\fm)$, and let $x_{1},\ldots,x_{d}$ be a system of parameters for $R$. Suppose that $T$ is an $R$-algebra, that $c$ is a non-zero divisor of $T$, while there is an $R$-linear map $\alpha:M \to T[c^{-1}]$. Let $M \to M'$ be a partial algebra modification of $M$ with respect to an initial segment of $x_{1},\ldots,x_{d}$, with degree bound $D$. Suppose that for every relation $x_{k+1}t_{k+1}=\sum_{i=1}^{k}x_{i}t_{i}$, $t_{i} \in T$, we have that $ct_{k+1} \in (x_{1},\ldots,x_{k})T$. Finally, suppose that $\alpha(M) \subseteq c^{-N}T$ for some integer $N>0$. Then the map $\alpha:M \to T[c^{-1}]$ fits into the commutative square:
$$
\begin{CD}
T[c^{-1}] @= T[c^{-1}] \\
@A\alpha AA @A\beta AA \\
M @>>> M'
\end{CD}
$$
in which $\beta:M' \to T[c^{-1}]$ is an $R$-linear map with image contained in $c^{-(ND+D+N)}T$.
\end{lemma}

\section{Almost purity theorem and almost Cohen-Macaulay algebras}

In this section, we establish the existence of a big Cohen-Macaulay algebra in the case that $R \to S$ is module-finite and $R[\frac{1}{p}] \to S[\frac{1}{p}]$ is \'etale and $R=\mathbf{W}(k)[[x_2,\ldots,x_d]]$ or $R=\mathbf{W}(k)[[t_1,\ldots,t_d]]/(p-G)$. We shall recall basic part of almost ring theory, for which we follow the treatment \cite{DavKed} and then make a definition of almost Cohen-Macaulay algebras. The advantage in working with almost ring theory by Davis and Kedlaya is that their approach is quite flexible and treats a general situation, so that it fits well into the search of big Cohen-Macaulay algebras in mixed characteristic. Some tricks to deal with the ramified case  may be found in the monograph \cite{GR2}. As \cite{GR2} is still undergoing revision, we give complete proofs of key facts for the convenience of readers. \cite{GR2} contains a detailed proof of some special cases of the Direct Summand Conjecture in the context of logarithmic geometry. We should emphasize that we consider the trivial log structure (the monoidal structure on $\mathcal{O}_X^{\times}$ for the local scheme $X=\Spec R$) in this article, but the powerful idea of the logarithmic geometry allows one to introduce some nice log structure (called "log regular") into a certain singular scheme, as if it were without any singular points. In general, if $(X, \mathcal{M}_X)$ is a regular log scheme, then the underlying scheme $X$ is known to be normal and Cohen-Macaulay (a proof of this fact may be found in \cite{GR2}).

\begin{definition}
\label{p-ideal}
Let $A$ be a $p$-torsion free ring. 

\begin{enumerate}
\item[$\mathrm{(i)}$]
A \textit{p-ideal} $I$ of $A$ is an ideal such that $I^m \subset p A$ for some $m>0$. 

\item[$\mathrm{(ii)}$]

An $A$-module $M$ is \textit{almost zero}, if $I M=0$ for any $p$-ideal $I$ of $A$.

\item[$\mathrm{(iii)}$]
An $A$-module $M$ is \textit{almost finite projective}, if for any given $p$-ideal $I \subset A$, there exists a finite free $A$-module $F$ for which there exists an $A$-module map $M \to F \to M$ given by multiplication by some $t \in A$ and such that $I \subset tA$.
\end{enumerate}
\end{definition}

For these definitions, we refer the reader to \cite{DavKed}. It will be necessary to prove that a certain algebra over a local ring $(R,\fm)$ is separated in the $\fm$-adic topology to construct a regular sequence out of an almost regular sequence.

\begin{lemma}
\label{separated}
Fix a prime integer $p>0$. Let $A$ be a $p$-torsion free ring and $A$ is $J$-adically separated for an ideal $J \subset A$. Assume that $B$ is a $p$-torsion free $A$-algebra. Consider the following conditions:

\begin{enumerate}
\item[$\mathrm{(i)}$]
$B$ is almost finite projective over $A$.

\item[$\mathrm{(ii)}$]
There is a composite map of $A$-modules: $
B \to F \to B$ such that $F$ is a finite free $A$-module and $B \to F \to B$ is multiplication by $p$. 

\item[$\mathrm{(iii)}$]
$B$ is $J$-adically separated.
\end{enumerate}
Then $\mathrm{(i)}$ implies $\mathrm{(ii)}$. $\mathrm{(ii)}$ implies $\mathrm{(iii)}$.
\end{lemma}

\begin{proof}
 $\mathrm{(i)}$ implies $\mathrm{(ii)}$: From the definition of almost finite projectivity, for a $p$-ideal $I=(p)$, there exists a map of $A$-modules: $B \to F \to B$ such that this is multiplication by some $t \in A$, $F$ is a finite free $A$-module and $I \subset tA$. Then we can write $p=at$ for some $a \in A$. By composing $B \to F \to B$ with the map $B \xrightarrow{a} B$, we get a map $B \to F \to B$ which is multiplication by $p$, as required.

$\mathrm{(ii)}$ implies $\mathrm{(iii)}$: 
Let $B \to F \to B$ be as in $\mathrm{(ii)}$ and consider the composite map:
$$
\bigcap_{n>0} J^n B \to \bigcap_{n>0} J^n F \to \bigcap_{n>0} J^n B
$$
which is multiplication by $p$. Since the $A$-module $F$ is finite free and $A$ is $J$-adically separated, we have $\bigcap_{n>0} J^n F=0$, which implies that $p \cdot \bigcap_{n>0} J^nB=0$. Since $B$ is $p$-torsion free, we get $\bigcap_{n>0} J^n B=0$ and $B$ is $J$-adically separated.
\end{proof}

We recall the definition of \textit{Witt-perfect} rings, due to Davis and Kedlaya.

\begin{definition}
Let $A$ be a commutative ring and fix a prime number $p$. Then $A$ is called \textit{Witt-perfect} (or \textit{p-Witt-perfect}), if the Witt-Frobenius map
$$
\mathbf{F}:\mathbf{W}_{p^n}(A) \to \mathbf{W}_{p^{n-1}}(A)
$$
is surjective for all $n \ge 2$.
\end{definition}

Note that the above definition makes sense for any ring $A$, regardless of the characteristic of $A$. Moreover, if $A$ has prime characteristic $p>0$, the Witt-Frobenius map $\mathbf{F}$ in the definition coincides with the lifting of the usual Frobenius $\Frob:A \to A$. A proof of this together with a concrete description of $\mathbf{F}$ in the general case is found in (\cite{DavKed}; Lemma 1.4). Since the above definition is not so intuitive, we recall the following more useful criterion.

\begin{lemma}
\label{criterion}
Let $B$ be a $p$-torsion free ring. Then the following conditions are equivalent:

\begin{enumerate}
\item[$\mathrm{(i)}$]
$B$ is a Witt-perfect ring.

\item[$\mathrm{(ii)}$]
$B$ satisfies the following properties:
\begin{enumerate}
\item[$\bullet$]
The Frobenius map on $B/pB$ is surjective.

\item[$\bullet$]
There exist $r,s \in B$ such that $r^p \equiv -p \mod psB$ and $s^N \in pB$ for some $N>0$.
\end{enumerate}
\end{enumerate}
\end{lemma}

\begin{proof}
This is (\cite{DavKed}; $(ii) \iff (xviii)+(xvii)$ of Corollary 3.3).
\end{proof}

We state the almost purity theorem (\cite{DavKed}; Theorem 5.2) based on the Witt-perfect condition. This theorem will play a crucial role later in the construction of almost Cohen-Macaulay algebras.

\begin{theorem}[Almost purity theorem]
\label{almostpure}
Let $B$ be a $p$-torsion free Witt-perfect ring which is integrally closed in $B[\frac{1}{p}]$. Let $B \to C$ be a ring homomorphism such that $B[\frac{1}{p}] \to C[\frac{1}{p}]$ is finite \'etale and let $C'$ be the integral closure of $B$ in $C[\frac{1}{p}]$. 

\begin{enumerate}
\item[$\mathrm{(i)}$]
The ring $C'$ is also Witt-perfect.

\item[$\mathrm{(ii)}$]
$C'$ is an almost finite projective $B$-module.
\end{enumerate}
\end{theorem}

We shall say that the ring map $B \to C$ satisfying the assumption of Theorem \ref{almostpure} is \textit{almost \'etale}. For a property concerning almost \'etale maps, we refer the reader to (\cite{Shim2}; Lemma 10.1). There is a simple criterion to be used to check when a given $p$-torsion free ring is Witt-perfect. We will use this criterion in the following discussion.

\begin{theorem}
\label{Witt-perfect}
Let $B$ be a $p$-torsion free algebra over a Witt-perfect ring $A$. Then $B$ is Witt-perfect if and only if the Frobenius map is surjective on $B/pB$.
\end{theorem}

\begin{proof}
Since $A$ is Witt-perfect by assumption, this is immediate from Lemma \ref{criterion}.
\end{proof}

\begin{definition}
A $p$-torsion free ring $B$ is called \textit{p-big}, if $B$ contains sequences of elements $\{\pi_n\}_{n \in \mathbf{N}}$  and $\{u_n\}_{n \in \mathbf{N}}$ such that $\pi_0=p$, every $u_n$ is a unit of $B$ and $\pi_{n+1}^p=\pi_nu_n$ for all $n \in \mathbf{N}$.
\end{definition}

For later use, we construct a certain valuation ring. Let $k$ be a perfect field of characteristic $p>0$ and choose a sequence $\{\pi_n\}_{n \in \mathbf{N}}$ as above in $\mathbf{W}(k)^+$ and $u_n=1$ for all $n \in \mathbf{N}$. We define
$$
V_{p^\infty}:=\bigcup_{n>0} \mathbf{W}(k)[\pi_n],
$$
In what follows, we may write $\pi_n$ as $p^{p^{-n}}$. The following lemma is essential for dealing with the unramified case.

\begin{lemma}
\label{valuation}
Let the notation be as above. Then $V_{p^\infty}$ is a Witt-perfect, non-discrete valuation ring of Krull dimension one.
\end{lemma}

\begin{proof}
From the construction, $V_{p^\infty}$ is the ascending union of discrete valuation rings $\mathbf{W}(k)[\pi_n]$ and its unique maximal ideal is not finitely generated. Thus, $V_{p^\infty}$ is a non-discrete valuation ring of Krull dimension one. Therefore, it remains to show that it is Witt-perfect. For this, we need to show that Lemma \ref{criterion} (ii) holds on $V_{p^\infty}$. By construction, the Frobenius map:
$$
\Frob:V_{p^\infty}/p V_{p^\infty} \to V_{p^\infty}/p V_{p^\infty}
$$
is surjective, so it only suffices to verify the second condition of Lemma \ref{criterion} (ii). First, assume $p>2$. Then the equality $(-\pi_1)^p=-p$ satisfies our demand. Next, assume $p=2$. In this case, we choose $r=\pi_1$, $s=2$ and $N=1$. This finishes the proof that $V_{p^{\infty}}$ is Witt-perfect.
\end{proof}

The valuation ring $V_{p^\infty}$ is called \textit{deeply ramified} in \cite{GR} and (non $p$-adically complete) \textit{perfectoid} in \cite{Sch}. We prove that $R_{p^{\infty}}$ has the desired properties.

\begin{proposition}
\label{WittCM}
$R_{p^{\infty}}$ is a Witt-perfect and balanced big Cohen-Macaulay $R$-algebra.
\end{proposition}

\begin{proof}
Since $R_{p^{\infty}}$ is faithfully flat over a regular local ring $R$, $R_{p^{\infty}}$ is a balanced big Cohen-Macaulay $R$-algebra. By Lemma \ref{WittNormal}, the Frobenius map is surjective on $R_{p^{\infty}}/pR_{p^{\infty}}$. To see that $R_{p^{\infty}}$ is Witt-perfect, it suffices to check the second condition of Lemma \ref{criterion} (ii). We need to consider the unramified and ramified cases separately.

First, assume that $p \notin \fm^2$. Then $R_{p^\infty}$ contains $V_{p^\infty}$ as a subring. Then apply Theorem \ref{Witt-perfect} to deduce that $R_{p^{\infty}}$ is Witt-perfect.

Next, assume that $p \in \fm^2$. We prove that there exist an element $\pi \in R_{p^{\infty}}$ and a unit $u \in R_{p^{\infty}}^{\times}$ such that
\begin{equation}
\label{need}
\pi^p=pu.
\end{equation}
Then we may write $p=\sum_{i=1}^n b_i b_i'$ with $b_i, b_i' \in \fm$. Let $\fm_{p^{\infty}}$ be the unique maximal ideal of $R_{p^{\infty}}$. Since the Frobenius map on $R_{p^{\infty}}/pR_{p^{\infty}}$ is surjective, we may find $d_i,d_i' \in R_{p^{\infty}}$ and $c_i,c_i' \in \fm_{p^{\infty}}$ such that
$$
b_i=c_i^p+pd_i,~b_i'={c_i'}^p+pd_i'.
$$
Substituting these equations into $p=\sum_{i=1}^n b_i b_i'$, we obtain $p(1-h)=\sum_{i=1}^n c_i^p {c_i'}^p$ with the property that $1-h \in R_{p^{\infty}}^{\times}$. Set
$$
f:=\sum_{i=1}^n c_i {c_i'}~\mbox{and}~g:=p^{-1}\Big(\sum_{i=1}^n c_i^p {c_i'}^p-(\sum_{i=1}^n c_i {c_i'})^p\Big).
$$
Since $c_i,  {c_i'} \in \fm_{p^{\infty}}$, we have $f,g \in \fm_{p^{\infty}}$. Hence $
p(1-h)=f^p+pg$ and $p(1-g-h)=f^p$. So setting $\pi:=f$ and $u:=1-g-h$ will suffice to obtain $(\ref{need})$. Since the Frobenius endomorphism is surjective on $R_{p^{\infty}}/pR_{p^{\infty}}$, we can find $v,w \in R_{p^{\infty}}$ such that $-u=v^p+pw$. Since $u$ is a unit and $pw$ is in the maximal ideal of $R_{p^{\infty}}$, the element $v^p$ must be a unit. Substituting the equation $-u=v^p+pw$ into $(\ref{need})$, we get
$$
\pi^p=-pv^p-p^2w,
$$
which gives
$$
\Big(\frac{\pi}{v}\Big)^p \equiv -p \pmod{p^2}.
$$
By letting
$$
s=p,~N=1~\mbox{and}~r=\frac{\pi}{v},
$$
the second condition of Lemma \ref{criterion} (ii) is satisfied and $R_{p^{\infty}}$ is Witt-perfect, as claimed.
\end{proof}

\begin{lemma}
The algebra $R_{p^{\infty}}$ is $p$-big. That is, it contains sequences $\{\pi_n\}$ and $\{u_n\}$ such that $\pi_0=p$, $u_n$ is a unit of $R_{p^{\infty}}$ and $\pi_{n+1}^p=\pi_n u_n$. In particular, $I_n:=(\pi_n)$ is a $p$-ideal.
\end{lemma}

\begin{proof}
When $R$ is unramified, this was already shown, as it contains a $p$-big valuation ring $V_{p^{\infty}}$. When $R$ is ramified, we may put $\pi_0:=p$ and $\pi_1:=\pi'$, where $\pi'$ is as in Proposition \ref{WittCM}. One may proceed further to define a sequence of elements $\{\pi_n\}$ with the required properties, proving that $R_{p^{\infty}}$ is $p$-big.
\end{proof}

The following definition of almost Cohen-Macaulay algebras is slightly different from the one given in \cite{Rob}. In fact, the $\fm$-adic separatedness is not assumed in \cite{Rob}. Our definition has the advantage that these conditions are sufficient to construct a big Cohen-Macaulay algebra in mixed characteristic.

\begin{definition}
Let $B$ be a $p$-torsion free $p$-big algebra over a Noetherian local ring $(R,\fm)$ of dimension $d>0$. Then $B$ is an \textit{almost Cohen-Macaulay R-algebra}, if the following conditions are satisfied:

\begin{enumerate}
\item[$\bullet$]
$B$ is separated in the $\fm$-adic topology.

\item[$\bullet$]
There is a system of parameters $x_1,\ldots,x_d$ of $R$ such that
$$
\frac{((x_1,\ldots,x_i):_Bx_{i+1})}{(x_1,\ldots,x_i)}
$$
is an almost zero $B$-module for all $i=0,\ldots,d-1$.
\end{enumerate}
\end{definition}

In particular, $\fm B \ne B$ from the definition.

\begin{proposition}
\label{algebramodification}
Suppose that $B$ is $p$-torsion free and an almost Cohen-Macaulay algebra over a local ring $(R,\fm)$. Assume that $C$ is a $B$-algebra such that $B \to C$ is torsion free and almost finite projective. Then the following assertions hold:

\begin{enumerate}
\item[$\mathrm{(i)}$]
$C$ is almost Cohen-Macaulay.

\item[$\mathrm{(ii)}$]
There exists a $C$-algebra $D$ such that $D$ is a big Cohen-Macaulay $R$-algebra.
\end{enumerate}
\end{proposition}

\begin{proof}
(i): The condition on almost Cohen-Macaulay property forces $B$ to be $\fm$-adically separated. Hence $C$ is $\fm$-adically separated in view of Lemma \ref{separated}. Moreover, $B$ is $p$-big. Fix a system of parameters $(x_1,\ldots,x_d)$ of $R$ which is an almost regular sequence on $B$. As our $p$-ideal, we choose $I_n=(\pi_n)$ and there exists a finite free $S$-module $F_n$ such that the composite map $C \to F_n \to C$ is multiplication by $t_n \in B$ with the property that $I_n \subset t_nB$ in view of the definition of almost projectivity. We may further take $t_n$ to generate a $p$-ideal by replacing it with some multiple. Since $F_n$ is $B$-free, it is almost Cohen-Macaulay. Applying the functor $(-)\otimes_C C/(x_1,\ldots,x_i)$, we get 
$$
C/(x_1,\ldots,x_i) \xrightarrow{f} F_n/(x_1,\ldots,x_i) \xrightarrow{g} C/(x_1,\ldots,x_i).
$$
Then since $\ker(f) \subset \ker(g \circ f)$ and $t_n \cdot \ker(g \circ f)=0$, we have $t_n \cdot \ker(f)=0$. Consider the commutative diagram:
$$
\begin{CD}
C/(x_1,\ldots,x_i) @>f>> F_n/(x_1,\ldots,x_i) \\
@VVx_{i+1}V @VVx_{i+1}V \\
C/(x_1,\ldots,x_i) @>f>> F_n/(x_1,\ldots,x_i)
\end{CD}
$$
Since $F_n$ is almost Cohen-Macaulay, we have 
$$
t_m \cdot \ker\big(F_n/(x_1,\ldots,x_i) \xrightarrow{x_{i+1}} F_n/(x_1,\ldots,x_i)\big)=0
$$
for any $m>0$. So the equality $t_n \cdot \ker(f)=0$, together with the snake lemma implies the following:
$$
t_m t_n \cdot ((x_1,\ldots,x_i):_Cx_{i+1}) \subset (x_1,\ldots,x_i).
$$
As both $m$ and $n$ may be taken arbitrarily large and each $t_n$ generates a $p$-ideal, we conclude that $C$ is almost Cohen-Macaulay.

(ii): We shall apply partial algebra modifications. We list the data of $C$ that is necessary to complete the proof.
\begin{enumerate}
\item[$\bullet$]
$C$ is $\fm$-adically separated and $p$-torsion free.

\item[$\bullet$]
$(x_1,x_2,\ldots,x_d)$ forms an almost regular sequence on $C$ in the sense defined just above.
\end{enumerate}
Then we may form the diagram consisting of sequences of partial algebra modifications with respect to a fixed system of parameters $(x_1,x_2,\ldots,x_d)$:
$$
\begin{CD}
C[\frac{1}{p}] @= C[\frac{1}{p}] @= \cdots @= C[\frac{1}{p}] \\
@AAA @AAA @. @AAA  \\
C @>>> M_1 @>>> \cdots @>>> M_r \\
\end{CD}
$$
by applying Lemma \ref{lemma3}. Assume that this sequence is bad. Then we may utilize the properties of $C$ stated as above and argue as in the proof of (\cite{Ho02}; Theorem 5.2) to derive a contradiction. Hence $C$ maps to a big Cohen-Macaulay algebra, as desired.
\end{proof}

We are now ready to prove the main theorem of this article.

\begin{theorem}
\label{maintheorem}
Let $R$ be a regular local ring of mixed characteristic $p>0$ and let $S$ be a torsion free module-finite $R$-algebra such that the localization $R[\frac{1}{p}] \to S[\frac{1}{p}]$ is \'etale. Then $S$ has a balanced big Cohen-Macaulay $R$-algebra.
\end{theorem}

\begin{proof}
First, there exists a flat local extension $(R,\fm) \to (R',\fm')$ such that $\fm R'=\fm'$ and the residue field $k$ of $R'$ is perfect. In particular, $\dim R=\dim R'$ and $R'$ is regular local. Now let $\widehat{R'}$ be the $\fm'$-adic completion of $R'$. Then $\widehat{R'}$ is isomorphic to $\mathbf{W}(k)[[x_2,\ldots,x_d]]$ or $\mathbf{W}(k)[[t_1,\ldots,t_d]]/(p-G)$. By \'etale base change, $\widehat{R'} \hookrightarrow \widehat{R'} \otimes_R S$ is \'etale after inverting $p$. Since $S \to \widehat{R'} \otimes_R S$ is a faithfully flat map between rings of the same Krull dimension, if $\widehat{R'} \otimes_R S$ maps to an algebra $\mathscr{B}(\widehat{R'} \otimes_R S)$ such that a system of parameters of $\widehat{R'}$ maps to a regular sequence in $\mathscr{B}(\widehat{R'} \otimes_R S)$, then $\mathscr{B}(\widehat{R'} \otimes_R S)$ is a big Cohen-Macaulay $R$-algebra. Then to prove the theorem, we may replace the original extension $R \hookrightarrow S$ with the extension: 
$$
\widehat{R'} \hookrightarrow (\widehat{R'} \otimes_R S).
$$
Recall that we have constructed Witt-perfect rings (see Lemma \ref{WittNormal} and Proposition \ref{WittCM} for their properties):
$$
R_{p^{\infty}}=\bigcup_{n>0} \mathbf{W}(k)[\pi_n][[x_2^{p^{-n}},\ldots,x_d^{p^{-n}}]]
$$
and
$$
R_{p^{\infty}}=\bigcup_{n>0} \mathbf{W}(k)[[t_1^{p^{-n}},\ldots,t_d^{p^{-n}}]]/(p-G)).
$$
Since $R[\frac{1}{p}] \to S[\frac{1}{p}]$ is \'etale, $R_{p^\infty}[\frac{1}{p}] \to (R_{p^\infty} \otimes_R S)[\frac{1}{p}]$ is also \'etale. In particular, $(R_{p^\infty} \otimes_R S)[\frac{1}{p}]$ is a normal ring which is a finite product of normal domains. We define 
$$
S_{p^\infty}=\mbox{the integral closure of}~R_{p^\infty}~\mbox{in}~(R_{p^\infty} \otimes_R S)[\frac{1}{p}].
$$ 
Then $S_{p^\infty}$ is a normal $S$-algebra. By applying both Theorem \ref{almostpure} and Proposition \ref{algebramodification}, we conclude that $S_{p^\infty}$ maps to a big Cohen-Macaulay $R$-algebra, whose $\fm$-adic completion is a balanced big Cohen-Macaulay $R$-algebra by (\cite{BH}; Corollary 8.5.3), which completes the proof of the theorem.
\end{proof}

As a corollary, we have the following.

\begin{corollary}
The normal $S$-algebra $S_{p^\infty}$ is almost Cohen-Macaulay.
\end{corollary}

We deduce another corollary.

\begin{corollary}
\label{maincorollary}
Let $R$ be a Noetherian regular domain and let $S$ be a torison free module-finite $R$-algebra. Assume that $R$ is $p$-torsion free and the localization $R[\frac{1}{p}] \to S[\frac{1}{p}]$ is finite \'etale for some prime integer $p>0$. Then $R \hookrightarrow S$ splits as an $R$-module homomorphism.
\end{corollary}

\begin{proof}
First, we observe that $R \to S$ splits if and only if $R_{\fp} \to S_{\fp}:=R_{\fp} \otimes_R S$ splits as an $R_{\fp}$-homomorphism for all prime ideals $\fp \in \Spec R$. If $p \notin \fp$, then it follows that $R_{\fp} \to S_{\fp}$ is \'etale. Therefore, this map splits for an obvious reason. So assume that $p \in \fp$. Then replacing the original extension $R \to S$ with $R_{\fp} \to S_{\fp}$, we may assume that $R$ is a regular local ring of mixed characteristic $p$. In particular, $S$ is a torsion free module-finite $R$-algebra and $S$ is semilocal. Then by Theorem \ref{maintheorem}, $S$ maps to a big Cohen-Macaulay $R$-algebra. Then by (\cite{Ho73}; Theorem 1), $R \to S$ splits if and only if we have $(x_1 \cdots x_d)^k \notin (x_1^{k+1},\ldots,x_d^{k+1})S$ for a regular system of parameters $x_1,\ldots,x_d$ of $R$ and every integer $k>0$. To prove the corollary by contradiction, assume that $(x_1 \cdots x_d)^k \in (x_1^{k+1} ,\ldots,x_d^{k+1})S$. Then mapping this relation to a big Cohen-Macaulay algebra yields a contradiction to the regularity of the sequence $x_1,\ldots,x_d$. Therefore, $R \to S$ must split, which completes the proof of the corollary.
\end{proof}

Finally, we mention that Bhatt proved that there is a complete local ring admitting no small Cohen-Macaulay algebra using Witt-vector cohomology in \cite{Bha2}. However, it is still an open question whether every complete local domain has a small Cohen-Macaulay module.

\begin{ack}
I would like to thank B. Bhatt, R. Heitmann, and P. Roberts for their careful reading. I also thank M. Asgharzadeh, Alberto F. Boix, and G. Piepmeyer for their comments. My special thanks go to L. E. Miller, who taught the author about Witt vectors. I am also grateful to an anonymous referee for finding errors and making suggestions.
\end{ack}


\begin{thebibliography}{99}



\bibitem{Ar} 
M. Artin, \emph{On the joins of Hensel rings}, Advances in Math. \textbf{7}  (1971),  282--296.



\bibitem{AsSh}
M. Asgharzadeh and K. Shimomoto, \emph{Almost Cohen-Macaulay and almost regular algebras via almost flat extensions}, J. Commutative Algebra \textbf{4} (2012), 445--478.


\bibitem{Bha1}
B. Bhatt,  \emph{Almost direct summands}, Nagoya Math. J. \textbf{214} (2014), 195--204.


\bibitem{Bha2}
B. Bhatt,  \emph{On the non-existence of small Cohen-Macaulay algebras}, J. Algebra \textbf{411} (2014), 1--11.


\bibitem{BH}
W. Bruns and J. Herzog,  \emph{Cohen-Macaulay rings}, Cambridge University Press, \textbf{39} (1998).


\bibitem{DavKed} 
C. Davis and K. S. Kedlaya, \emph{On the Witt vector Frobenius}, Proceedings of the A.M.S. \textbf{142} (2014), 2211--2226.


\bibitem{GR}
O. Gabber and L. Ramero, \emph{Almost ring theory}, LNM \textbf{1800}, Springer-Verlag  (2003).


\bibitem{GR2}
O. Gabber and L. Ramero, \emph{Foundations for almost ring theory}, \textbf{arXiv:math/0409584}.


\bibitem{He1}
R. Heitmann, \emph{The direct summand conjecture in dimension three}, Ann. of Math. \textbf{156} (2002), 695--712.


\bibitem{Ho73}
M. Hochster,  \emph{Contracted ideals from integral extensions of regular rings}, Nagoya Math. J. \textbf{51} (1973), 25--43.


\bibitem{Ho02}
M. Hochster,  \emph{Big Cohen-Macaulay algebras in dimension three via Heitmann's theorem}, J. Algebra \textbf{254} (2002), 395--408.


\bibitem{Ho07}
M. Hochster, \emph{Homological conjectures, old and new},  Illinois J. Math. \textbf{51} (2007), 151--169. 



\bibitem{KedRuo}
K. S. Kedlaya and R. Liu, \emph{Relative p-adic Hodge theory, I: Foundations}, to appear in Ast\'erisque.


\bibitem{Mat}
H. Matsumura, \emph{Commutative ring theory}, Cambridge University
Press. \textbf{8} (1986).


\bibitem{Oli}
J. P. Olivier, \emph{Going up along absolutely flat morphisms}, J. Pure and Applied Algebra \textbf{30} (1983), 47--59.


\bibitem{Rab}
J. Rabinoff, \emph{The theory of Witt vectors}, \textbf{arXiv:1409.7445}.


\bibitem{Rob}
P. Roberts, \emph{Fontaine rings and local cohomology}, J. Algebra \textbf{323} (2010), 2257--2269.


\bibitem{Rob12}
P. Roberts, \emph{The homological conjectures}, Progress in Commutative Algebra 1: Combinatorics and Homology (de Gruyter Proceedings in Mathematics)  Sean Sather-Wagstaff, Christopher Francisco, Lee C. Klingler (2012).


\bibitem{Sch}
P. Scholze,  \emph{Perfectoid spaces}, Publ. Math. de l'IH\'ES \textbf{116} (2012), 245--313.


\bibitem{Se}
J.-P. Serre,  \emph{Local Fields}, Graduate Texts in Mathematics \textbf{67}, Springer-Verlag, New York (1973).


\bibitem{Shim}
K. Shimomoto, \emph{Almost Cohen-Macaulay algebras in mixed characteristic via Fontaine rings}, Illinois J. Math. \textbf{55} (2011), 107--125.

\bibitem{Shim2}
K. Shimomoto, \emph{An embedding problem of Noetherian rings into the Witt vectors}, preprint.


\bibitem{SwHu}
I. Swanson and C. Huneke, \emph{Integral closure of ideals rings, and modules}, London Math. Society Lecture Note Series \textbf{336}.
\end{thebibliography}
\end{document}